\documentclass[11pt,draft]{article}
\usepackage{amscd}
\usepackage{amsmath}
\usepackage{latexsym}
\usepackage{amsfonts}
\usepackage{amssymb}
\usepackage{mathrsfs}
\usepackage{amsthm}

\newtheorem{theorem}{Theorem}[section]
\newtheorem{lemma}{Lemma}[section]

\newtheorem{proposition}{Proposition}[section]
\newtheorem{definition}{Definition}[section]

\newtheorem{remark}{Remark}[section]

\numberwithin{equation}{section}

\DeclareMathOperator{\intd}{d}
\DeclareMathOperator{\Lip}{Lip}
\DeclareMathOperator{\Div}{div}
\begin{document}

\title{Global well-posedness of the critical Burgers equation in critical Besov spaces}

\author{Changxing Miao$^1$ and  Gang Wu$^2$\\
        \\
        \small{$^{1}$ Institute of Applied Physics and Computational Mathematics,}\\
        \small{P.O. Box 8009, Beijing 100088, P.R. China.}\\
        \small{(miao\_{}changxing@iapcm.ac.cn)}\\
        \small{$^2$  The Graduate School of China Academy of Engineering Physics,}\\
        {\small P.O. Box 2101, Beijing 100088, P.R. China.}\\
        {\small{(wugangmaths@yahoo.com.cn)}}}
\date{}
\maketitle

\begin{abstract}
We make use of the method of modulus of continuity \cite{K-N-S}
and Fourier localization technique \cite{A-H} to prove the global
well-posedness of the critical Burgers equation
$\partial_{t}u+u\partial_{x}u+\Lambda u=0$ in critical Besov
spaces $\dot{B}^{\frac{1}{p}}_{p,1}(\mathbb{R})$ with
$p\in[1,\infty)$, where $\Lambda=\sqrt{-\triangle}$.
\end{abstract}

2000 Mathematics Subject Classification: 35K55, 35Q53

Key words and phrases: Burgers equation; Modulus of continuity; Fourier localization; Global well-posedness; Besov spaces
\section{Introduction}
We consider the Burgers equation with fractional dissipation in $\mathbb{R}$,
\begin{equation}\label{eq1.1}
    \begin{cases}
        \partial_{t}u+u\partial_{x}u+\Lambda^{\alpha}u=0\\
        u(x,0)=u_{0}(x),
    \end{cases}
\end{equation}
where $0\leq\alpha\leq 2$ and the operator $\Lambda^{\alpha}$ is defined by Fourier transform
\begin{equation*}
    \mathcal{F}(\Lambda^{\alpha}u)(\xi)=|\xi|^{\alpha}\mathcal{F}u(\xi).
\end{equation*}

The Burgers equation \eqref{eq1.1} with $\alpha=0$ and $\alpha=2$ has received
an extensive amount of attention since the studies by Burgers in
the 1940s. If $\alpha=0$, the equation is perhaps the most basic
example of a PDE evolution leading to shocks; if $\alpha=2$, it
provides an accessible model for studying the interaction between
nonlinear and dissipative phenomena. Recently, in \cite{K-N-S} for
the periodic case authors give a complete study for general $\alpha\in[0,2]$, see also \cite{A-D,Dong2,K-M-X,M-Y-Z}.
In particular, for $\alpha=1$, with help of the method of modulus of
continuity they proved the global well-posedness of the equation
in the critical Hilbert space $H^{\frac{1}{2}}(\mathbb{T}^{1})$.

In this paper, we study the following critical case,
\begin{equation}\label{eq1.2}
    \begin{cases}
        \partial_{t}u+u\partial_{x}u+\Lambda u=0\\
        u(x,0)=u_{0}(x).
    \end{cases}
\end{equation}
We use similar arguments as in \cite{A-H}.
Making use of Fourier localization technique and the method of modulus of continuity \cite{K-N-S},
we prove the global well-posedness of the critical Burgers equation \eqref{eq1.2}
in critical Besov spaces $\dot{B}^{\frac{1}{p}}_{p,1}(\mathbb{R})$ with $p\in[1,\infty)$.

It is well known that $\dot{B}^{\frac{1}{p}}_{p,1}$ is the critical space under the scaling invariance.
That is, if $u(x,t)$ is a solution of \eqref{eq1.2},
then $u_{\lambda}(x,t)=u(\lambda x,\lambda t)$ is also a solution of the same equation
and $\|u_{\lambda}(\cdot,t)\|_{\dot{B}^{\frac{1}{p}}_{p,1}}\approx\|u(\cdot,\lambda t)\|_{\dot{B}^{\frac{1}{p}}_{p,1}}$.

Now we give out our main results. The first main result is the
following:
\begin{theorem}\label{thm1.1}
Let $u_{0}\in \dot{B}^{\frac{1}{p}}_{p,1}(\mathbb{R})$ with $p\in[1,\infty)$,
then the critical Burgers equation \eqref{eq1.2} has a unique global solution $u$ such that
\begin{equation*}
    u\in \mathcal{C}(\mathbb{R}^{+};\dot{B}^{\frac{1}{p}}_{p,1})
    \cap L^{1}_{loc}(\mathbb{R}^{+};\dot{B}^{\frac{1}{p}+1}_{p,1}).
\end{equation*}
\end{theorem}

\begin{remark}\label{rem1.1}
Because of the restriction of the smooth index $s$ stemming from the \emph{a priori} estimate
for the transport-diffusion equation (see Theorem \ref{thm1.2}),
we can not get the result for the limit case $p=\infty$.
\end{remark}

\begin{remark}\label{rem1.2}\;
The corresponding question for the quasi-geostrophic equation
has been a focus of significant effort (see e.g. \cite{A-H,C-M-Z,C-C-W,Dong1,K-N-V,Wu1,Wu2})
and the critical Q-G equation has been recently resolved in \cite{K-N-V} for periodic case.
Based on \cite{K-N-V}, Abidi-Hmidi in \cite{A-H} and Dong-Du in \cite{Dong1}
give the corresponding result for Cauchy problem of the critical Q-G equation in the  framework of Besov space
and Sobolev space, respectively.
After the present paper is completed, Prof. J.Wu and H.Dong informed us that the authors in \cite{Dong2}
gave the global well-posedness for the critical fractal Burgers equation in inhomogeneous space $H^{\frac12}$
by similar argument in \cite{Dong1}.
\end{remark}

In order to prove this theorem, we first prove the local well-posedness which is the major part of this paper.
Next we make use of the modulus of continuity \cite{K-N-S} to get the global well-posedness.
We mention that the property allowing us to remove the periodicity is the spatial decay of the solution.

The key of proving the local well-posedness is an optimal \emph{a priori} estimate
for the following transport-diffusion equation in $\mathbb{R}^{N}$:
\begin{equation*}
    \begin{cases}
        \partial_{t}u+v\cdot\nabla u+\nu\Lambda^{\alpha}u=f\\
        u(x,0)=u_{0}(x),
    \end{cases}
    \eqno{(TD)_{\nu,\alpha}}
\end{equation*}
where $v$ is a given vector field which needs not to be divergence free, $u_{0}$ is the initial data,
$f$ is a given external force term,
$\nu\geq 0$ is a constant, $0\leq\alpha\leq 2$.
Our second main result is the following:
\begin{theorem}\label{thm1.2}
Let $1\leq \rho_{1}\leq \rho\leq \infty$, $1\leq p\leq p_{1}\leq \infty$ and $1\leq r\leq \infty$.
Let $s\in \mathbb{R}$ satisfy the following
\begin{align*}
    & s<1+\frac{N}{p_{1}}\; \Big(\text{or}\; s\leq 1+\frac{N}{p_{1}}\; \text{if} \; r=1\Big),\\
    & s>-N\min\Big(\frac{1}{p_{1}},\frac{1}{p'}\Big)\;
    \bigg(\text{or}\; s>-1-N\min\Big(\frac{1}{p_{1}},\frac{1}{p'}\Big)\; \text{if}\; \Div v=0\bigg).
\end{align*}
There exists a constant $C>0$ depending only on $N$, $\alpha$, $s$, $p$, $p_{1}$ and $r$,
such that for any smooth solution $u$ of $(TD)_{\nu,\alpha}$ with $\nu\geq 0$, we have the following a priori estimate:
\begin{equation}\label{eq1.3}
    \nu^{\frac{1}{\rho}}\|u\|_{\widetilde{L}^{\rho}_{T}\dot{B}^{s+\frac{\alpha}{\rho}}_{p,r}}
    \leq Ce^{CZ(T)}\Big(\|u_{0}\|_{\dot{B}^{s}_{p,r}}+\nu^{\frac{1}{\rho_{1}}-1}
    \|f\|_{\widetilde{L}^{\rho_{1}}_{T}\dot{B}^{s-\alpha+\frac{\alpha}{\rho_{1}}}_{p,r}}\Big)
\end{equation}
with $Z(T):=\int_{0}^{T}\|\nabla v(t)\|_{\dot{B}^{\frac{N}{p_{1}}}_{p_{1},\infty}\cap L^{\infty}}\intd t$.

Besides if $u=v$, then for all $s>0$($s>-1$ if $\Div v=0$), the estimate \eqref{eq1.3} holds with
$Z(T)=\int_{0}^{T}\|\nabla v(t)\|_{L^{\infty}}\intd t$.
\end{theorem}

\begin{remark}\label{rem1.3}
When $\alpha=2$, the above a priori estimate has been proved by
R.Danchin in \cite{Danchin}. In this paper, we extend Danchin's
results to the general case $\alpha\in [0,2]$. The proof's key is
the use of Lagrangian coordinates transformation together with an
important commutator estimate.
\end{remark}

The rest of this paper is arranged as follows:

In Section 2, we recall some definitions and properties about homogeneous Besov spaces,
and we will also list some useful lemmas.
In Section 3, we prove Theorem \ref{thm1.2}.
In Section 4, we prove the local well-posedness.
In Section 5, we give the blow-up criterion.
In Section 6, we complete the proof of the global well-posedness.

\paragraph{Notation:}
Throughout the paper, $C$ stands for a constant which may be different in each occurrence.
We shall sometimes use the notation $A\lesssim B$ instead of $A\leq CB$ and $A\approx B$ means that $A\lesssim B$
and $B\lesssim A$.
\section{Preliminaries}
Let us first recall the Littlewood-Paley Theory.
Let $\chi$ and $\varphi$ be a couple of smooth radial functions valued in $[0,1]$
such that $\chi$ is supported in the ball $\big\{\xi\in\mathbb{R}^{N}\big||\xi|\leq\frac{4}{3}\big\}$,
$\varphi$ is supported in the shell $\big\{\xi\in\mathbb{R}^{N}\big|\frac{3}{4}\leq|\xi|\leq\frac{8}{3}\big\}$
and
\begin{align*}
    \chi(\xi)+\sum_{q\in \mathbb{N}}\varphi(2^{-q}\xi)=1,\quad \forall\xi\in \mathbb{R}^{N};\\
    \sum_{q\in \mathbb{Z}}\varphi(2^{-q}\xi)=1,\quad \forall\xi\in \mathbb{R}^{N}\backslash\{0\}.
\end{align*}
Denoting $\varphi_{q}(\xi)=\varphi(2^{-q}\xi)$ and $h_{q}=\mathcal{F}^{-1}\varphi_{q}$,
we define the homogeneous dyadic blocks as
\begin{equation*}
    \dot{\Delta}_{q}u:=\varphi(2^{-q}D)u=\int_{\mathbb{R}^{N}}h_{q}(y)u(x-y)\intd y,\quad \forall q\in \mathbb{Z}.
\end{equation*}
We can also define the following low-frequency cut-off:
\begin{equation*}
    \dot{S}_{q}u:=\sum_{j\leq q-1}\dot{\Delta}_{j}u.
\end{equation*}

\begin{definition}\label{def2.1}
Let $\mathcal{S}_{h}'$ be the space of temperate distributions $u$ such that
\begin{equation*}
    \lim_{q\rightarrow-\infty}\dot{S}_{q}u=0,\quad \text{in} \quad \mathcal{S}'.
\end{equation*}
\end{definition}

The formal equality
\begin{equation*}
    u=\sum_{q\in \mathbb{Z}}\dot{\Delta}_{q}u
\end{equation*}
holds in $\mathcal{S}_{h}'$ and is called the \emph{homogeneous Littlewood-Paley decomposition}.
It has nice properties of quasi-orthogonality:
\begin{equation}\label{eq2.1}
    \dot{\Delta}_{q'}\dot{\Delta}_{q}u\equiv 0\quad \text{if}\quad |q'-q|\geq 2\quad\text{and}\quad
    \dot{\Delta}_{q'}(\dot{S}_{q-1}u\dot{\Delta}_{q}v)\equiv 0\quad \text{if}\quad |q'-q|\geq 5.
\end{equation}

Let us now define the homogeneous Besov spaces:
\begin{definition}\label{def2.2}
For $s\in \mathbb{R}$, $(p,r)\in [1,\infty]^{2}$ and $u\in \mathcal{S}_{h}'$, we set
\begin{equation*}
    \|u\|_{\dot{B}^{s}_{p,r}}
    :=\Big(\sum_{q\in \mathbb{Z}}2^{qsr}\|\dot{\Delta}_{q}u\|_{L^{p}}^{r}\Big)^{\frac{1}{r}}
    \quad\text{if}\quad r<\infty
\end{equation*}
and
\begin{equation*}
    \|u\|_{\dot{B}^{s}_{p,\infty}}:=\sup_{q\in\mathbb{Z}}2^{qs}\|\dot{\Delta}_{q}u\|_{L^{p}}.
\end{equation*}
We then define the \emph{homogeneous Besov spaces} as
\begin{equation*}
    \dot{B}^{s}_{p,r}:=\big\{u\in\mathcal{S}_{h}'\big|\|u\|_{\dot{B}^{s}_{p,r}}<\infty\big\}.
\end{equation*}
\end{definition}

The above definition does not depend on the choice of the couple $(\chi,\varphi)$.
Remark that if $s<\frac{N}{p}$ or $s=\frac{N}{p}$ and $r=1$, then $\dot{B}^{s}_{p,r}$ is a Banach space.

We now recall some basic properties of the homogeneous Besov spaces.
\begin{proposition}\label{prop2.1}
The following properties hold true(cf. \cite{R-S,Triebel}):
\begin{enumerate}
\item Generalized derivatives: Let $\sigma\in\mathbb{R}$, then the operator $\Lambda^{\sigma}$ is an isomorphism
    from $\dot{B}^{s}_{p,r}$ to $\dot{B}^{s-\sigma}_{p,r}$.
\item Sobolev embedding: If $p_{1}\leq p_{2}$ and $r_{1}\leq r_{2}$,
    then $\dot{B}^{s}_{p_{1},r_{1}}\hookrightarrow \dot{B}^{s-N(\frac{1}{p_{1}}-\frac{1}{p_{2}})}_{p_{2},r_{2}}$.
\item If $(p,r)\in [1,\infty]^{2}$ and $s>0$, there exists a positive constant $C=C(N,s)$ such that
    \begin{equation*}
        \|uv\|_{\dot{B}^{s}_{p,r}}\leq C(\|u\|_{L^{\infty}}\|v\|_{\dot{B}^{s}_{p,r}}
            +\|v\|_{L^{\infty}}\|u\|_{\dot{B}^{s}_{p,r}}).
    \end{equation*}
\end{enumerate}
\end{proposition}

In our next study we require two kinds of coupled space-time Besov spaces.
The first one is defined by the following manner: for $T>0$ and $\rho\in [1,\infty]$,
we denote by $L^{\rho}_{T}\dot{B}^{s}_{p,r}$ the set of all tempered distribution $u$ satisfying
\begin{equation*}
    \|u\|_{L^{\rho}_{T}\dot{B}^{s}_{p,r}}:=
    \bigg\|\Big(\sum_{q\in \mathbb{Z}}2^{qsr}\|\dot{\Delta}_{q}u\|_{L^{p}}^{r}\Big)^{\frac{1}{r}}\bigg\|
    _{L^{\rho}_{T}}<\infty.
\end{equation*}
The second mixed space is $\widetilde{L}^{\rho}_{T}\dot{B}^{s}_{p,r}$
which is the set of all tempered distribution $u$ satisfying
\begin{equation*}
    \|u\|_{\widetilde{L}^{\rho}_{T}\dot{B}^{s}_{p,r}}
    :=\Big(\sum_{q\in \mathbb{Z}}2^{qsr}\|\dot{\Delta}_{q}u\|_{L^{\rho}_{T}L^{p}}^{r}\Big)^{\frac{1}{r}}<\infty.
\end{equation*}
Let us remark that, by virtue of the Minkowski inequality, we have
\begin{equation*}
    \|u\|_{\widetilde{L}^{\rho}_{T}\dot{B}^{s}_{p,r}}\leq \|u\|_{L^{\rho}_{T}\dot{B}^{s}_{p,r}}
    \quad \text{if}\quad \rho\leq r,
\end{equation*}
and
\begin{equation*}
    \|u\|_{L^{\rho}_{T}\dot{B}^{s}_{p,r}}\leq \|u\|_{\widetilde{L}^{\rho}_{T}\dot{B}^{s}_{p,r}}
    \quad \text{if}\quad \rho\geq r.
\end{equation*}

Now we give some useful lemmas.
\begin{lemma}\label{lem2.2}(cf. \cite{H-K,W-Y})
Let $\phi$ be a smooth function supported in the shell
$\big\{\xi\in\mathbb{R}^{N}\big|R_{1}\leq |\xi|\leq R_{2}, 0<R_{1}<R_{2}\big\}$. There exist two positive constants
$\kappa$ and $C$ depending only on $\phi$ such that for all $1\leq p\leq \infty$, $\tau\geq 0$ and $\lambda>0$,
we have
\begin{equation*}
    \|\phi(\lambda^{-1}D)e^{-\tau\Lambda^{\alpha}}u\|_{L^{p}}
    \leq Ce^{-\kappa\tau\lambda^{\alpha}}\|\phi(\lambda^{-1}D)u\|_{L^{p}}.
\end{equation*}
\end{lemma}

\begin{lemma}\label{lem2.3}(cf. \cite{Danchin})
Let $v$ be a smooth vector field. Let $\psi_{t}$ be the solution to
\begin{equation*}
    \psi_{t}(x)=x+\int_{0}^{t}v(\tau,\psi_{\tau}(x))\intd \tau.
\end{equation*}
Then for all $t\in\mathbb{R}^{+}$, the flow $\psi_{t}$ is a $C^{1}$ diffeomorphism over $\mathbb{R}^{N}$ and one has
\begin{align*}
    &\|\nabla\psi_{t}^{\pm 1}\|_{L^{\infty}}\leq e^{V(t)},\\
    &\|\nabla\psi_{t}^{\pm 1}-{\rm Id}\|_{L^{\infty}}\leq e^{V(t)}-1,\\
    &\|\nabla^{2}\psi_{t}^{\pm 1}\|_{L^{\infty}}\leq e^{V(t)}\int_{0}^{t}\|\nabla^{2}v(\tau)\|_{L^{\infty}}
        e^{V(\tau)}\intd \tau,
\end{align*}
where $V(t)=\int_{0}^{t}\|\nabla v(\tau)\|_{L^{\infty}}\intd \tau$.
\end{lemma}

\begin{lemma}\label{lem2.4}(cf. \cite{C-M-W})
Let $v$ be a given vector field belonging to $L^{1}_{loc}(\mathbb{R}^{+};\Lip)$.
For $q\in\mathbb{Z}$ we set $u_{q}:=\dot{\Delta}_{q}u$
and denote by $\psi_{q}$ the flow of the regularized vector field $\dot{S}_{q-1}v$.
Then for $u\in\dot{B}^{\alpha}_{p,\infty}$ with $\alpha\in[0,2)$ and $p\in[1,\infty]$ we have
\begin{equation*}
    \|\Lambda^{\alpha}(u_{q}\circ\psi_{q})-(\Lambda^{\alpha}u_{q})\circ\psi_{q}\|_{L^{p}}
    \leq Ce^{CV(t)}V^{1-\frac{\alpha}{2}}(t)2^{q\alpha}\|u_{q}\|_{L^{p}},
\end{equation*}
where $V(t)=\int_{0}^{t}\|\nabla v(\tau)\|_{L^{\infty}}\intd \tau$ and $C=C(\alpha,p)>0$ is a constant.
\end{lemma}

\begin{lemma}\label{lem2.5}
Let $\sigma\in\mathbb{R}$ and $1\leq p\leq p_{1}\leq\infty$, $p_{2}:=(1/p-1/p_{1})^{-1}$.
Let $R_{q}:=(\dot{S}_{q-1}v-v)\cdot\nabla \dot{\Delta}_{q}u-[\dot{\Delta}_{q},v\cdot\nabla]u$.
There exists a constant $C=C(N,\sigma)$ such that
\begin{equation*}
    \begin{split}
      2^{q\sigma}\|R_{q}\|_{L^{p}} & \leq C\bigg(\sum_{|q'-q|\leq 4}\|\dot{S}_{q'-1}\nabla v\|_{L^{\infty}}
            2^{q'\sigma}\|\dot{\Delta}_{q'}u\|_{L^{p}}\\
        & \quad+\sum_{q'\geq q-3}2^{q-q'}\|\dot{\Delta}_{q'}\nabla v\|_{L^{\infty}}
            2^{q\sigma}\|\dot{\Delta}_{q}u\|_{L^{p}} \\
        & \quad+\sum_{\substack{|q'-q|\leq 4\\q''\leq q'-2}}2^{(q-q'')(\sigma-1-\frac{N}{p_{1}})}
            2^{q'\frac{N}{p_{1}}}\|\dot{\Delta}_{q'}\nabla v\|_{L^{p_{1}}}
            2^{q''\sigma}\|\dot{\Delta}_{q''}u\|_{L^{p}} \\
        & \quad+\sum_{\substack{q'\geq q-3\\|q''-q'|\leq 1}}
            2^{(q-q')\big(\sigma+N\min(\frac{1}{p_{1}},\frac{1}{p'})\big)}\\
        & \qquad \times 2^{q'\frac{N}{p_{1}}}
            \big(2^{q-q'}\|\dot{\Delta}_{q'}\nabla v\|_{L^{p_{1}}}+\|\dot{\Delta}_{q'}\Div v\|_{L^{p_{1}}}\big)
            2^{q''\sigma}\|\dot{\Delta}_{q''}u\|_{L^{p}}\bigg),
    \end{split}
\end{equation*}
and the third term in the right-hand side may be replaced by
\begin{equation*}
      C\sum_{|q'-q|\leq 4}2^{q'(\sigma-1)}\|\dot{\Delta}_{q'}\nabla v\|_{L^{p_{1}}}
            \|\dot{S}_{q'-1}\partial_{j}u\|_{L^{p_{2}}}.
\end{equation*}
Besides if $u=v$, the following estimate holds true:
\begin{equation*}
    \begin{split}
      2^{q\sigma}\|R_{q}\|_{L^{p}} & \leq C\bigg(\sum_{|q'-q|\leq 4}\|\dot{S}_{q'-1}\nabla u\|_{L^{\infty}}
            2^{q'\sigma}\|\dot{\Delta}_{q'}u\|_{L^{p}} \\
        & \quad+\sum_{q'\geq q-3}2^{q-q'}\|\dot{\Delta}_{q'}\nabla u\|_{L^{\infty}}
            2^{q\sigma}\|\dot{\Delta}_{q}u\|_{L^{p}} \\
        & \quad+\sum_{\substack{q'\geq q-3\\|q''-q'|\leq 1}}2^{(q-q')\sigma}\\
        & \qquad \times \big(2^{q-q'}\|\dot{\Delta}_{q'}\nabla u\|_{L^{\infty}}
            +\|\dot{\Delta}_{q'}\Div u\|_{L^{\infty}}\big)2^{q''\sigma}\|\dot{\Delta}_{q''}u\|_{L^{p}}\bigg).
    \end{split}
\end{equation*}
\end{lemma}
R.Dancin in \cite{Danchin1} gave the proof for the nonhomogeneous case.
For the convenience of the reader, we will give the proof for the homogeneous case in the appendix.
\section{Proof of Theorem \ref{thm1.2}}
\begin{proof}[Proof of Theorem \ref{thm1.2}]
Here we only prove the case $\alpha\in[0,2)$ (for the case $\alpha=2$, see \cite{Danchin}).

Let $u_{q}:=\dot{\Delta}_{q}u$ and $f_{q}:=\dot{\Delta}_{q}f$.
Applying $\dot{\Delta}_{q}$ to $(TD)_{\nu,\alpha}$ yields
\begin{equation*}
    \partial_{t}u_{q}+\dot{S}_{q-1}v\cdot\nabla u_{q}+\nu\Lambda^{\alpha}u_{q}=f_{q}+R_{q}
\end{equation*}
with $R_{q}:=(\dot{S}_{q-1}v-v)\cdot\nabla u_{q}-[\dot{\Delta}_{q},v\cdot\nabla]u$.

Let $\psi_{q}$ be the flow of the regularized vector field $\dot{S}_{q-1}v$.
Denote $\bar{u}_{q}:=u_{q}\circ\psi_{q}$, $\bar{f}_{q}:=f_{q}\circ\psi_{q}$ and $\bar{R}_{q}:=R_{q}\circ\psi_{q}$.
Then we have
\begin{equation}\label{eq3.1}
    \partial_{t}\bar{u}_{q}+\nu\Lambda^{\alpha}\bar{u}_{q}=\bar{f}_{q}+\bar{R}_{q}+\nu G_{q}
\end{equation}
with $G_{q}:=\Lambda^{\alpha}(u_{q}\circ\psi_{q})-(\Lambda^{\alpha}u_{q})\circ\psi_{q}$.

Applying $\dot{\Delta}_{j}$ to \eqref{eq3.1} and using Lemma \ref{lem2.2}, we get
\begin{equation}\label{eq3.2}
    \begin{split}
      \|\dot{\Delta}_{j}\bar{u}_{q}(t)\|_{L^{p}} & \lesssim e^{-\kappa\nu t2^{j\alpha}}
        \|\dot{\Delta}_{j}u_{0,q}\|_{L^{p}}+\int_{0}^{t}e^{-\kappa\nu (t-\tau)2^{j\alpha}}\\
        &\quad\times\big(\|\dot{\Delta}_{j}\bar{f}_{q}\|_{L^{p}}+\|\dot{\Delta}_{j}\bar{R}_{q}\|_{L^{p}}
        +\nu\|\dot{\Delta}_{j}G_{q}\|_{L^{p}}\big)\intd\tau.
    \end{split}
\end{equation}
Now from Lemma \ref{lem2.4} we have
\begin{equation}\label{eq3.3}
    \|\dot{\Delta}_{j}G_{q}(t)\|_{L^{p}}\leq Ce^{CV(t)}V^{1-\frac{\alpha}{2}}(t)2^{q\alpha}\|u_{q}\|_{L^{p}}.
\end{equation}
According to Bernstein lemma and Lemma \ref{lem2.3}, we can get
\begin{equation}\label{eq3.4}
    \begin{split}
    \|\dot{\Delta}_{j}\bar{f}_{q}(t)\|_{L^{p}} & \lesssim 2^{-j}\|\nabla\dot{\Delta}_{j}\bar{f}_{q}\|_{L^{p}}\\
    & \lesssim 2^{-j}\|(\nabla f_{q})\circ\psi_{q}\|_{L^{p}}\|\nabla\psi_{q}\|_{L^{\infty}}\\
    & \lesssim 2^{-j}\|\nabla f_{q}\|_{L^{p}}\|J_{\psi_{q}^{-1}}\|_{L^{\infty}}^{\frac{1}{p}}
        \|\nabla\psi_{q}\|_{L^{\infty}}\\
    & \lesssim e^{CV(t)}2^{q-j}\|f_{q}\|_{L^{p}}.
    \end{split}
\end{equation}
Arguing similarly as in deriving \eqref{eq3.4}, we obtain
\begin{equation*}
    \|\dot{\Delta}_{j}\bar{R}_{q}(t)\|_{L^{p}}\lesssim e^{CV(t)}2^{q-j}\|R_{q}\|_{L^{p}}.
\end{equation*}
According to Lemma \ref{lem2.5}, we get
\begin{equation}\label{eq3.5}
    \|\dot{\Delta}_{j}\bar{R}_{q}(t)\|_{L^{p}}\lesssim e^{CV(t)}2^{q-j}c_{q}(t)2^{-qs}Z'(t)\|u(t)\|_{\dot{B}^{s}_{p,r}},
\end{equation}
with $\|c_{q}(t)\|_{\ell^{r}}=1$.

Plugging \eqref{eq3.3}, \eqref{eq3.4} and \eqref{eq3.5} into \eqref{eq3.2}, taking the $L^{\rho}$ norm over $[0,t]$
and multiplying both sides by $\nu^{\frac{1}{\rho}}2^{q(s+\frac{\alpha}{\rho})}$, we obtain
\begin{equation}\label{eq3.6}
    \begin{split}
      \nu^{\frac{1}{\rho}}2^{q(s+\frac{\alpha}{\rho})}\|\dot{\Delta}_{j}\bar{u}_{q}\|_{L^{\rho}_{t}L^{p}} & \lesssim
            2^{\frac{(q-j)\alpha}{\rho}}(1-e^{-\kappa\nu \rho t2^{j\alpha}})^{\frac{1}{\rho}}
            2^{qs}\|\dot{\Delta}_{j}u_{0,q}\|_{L^{p}}\\
        & \quad +\nu^{-\frac{1}{\rho_{1}'}}2^{(q-j)(1+\frac{\alpha}{\rho}+\frac{\alpha}{\rho_{1}'})}
            e^{CV(t)}2^{q(s-\frac{\alpha}{\rho_{1}'})}\|f_{q}\|_{L^{\rho_{1}}_{t}L^{p}}\\
        & \quad +\nu^{\frac{1}{\rho}}2^{q(s+\frac{\alpha}{\rho})}2^{(q-j)\alpha}
            e^{CV(t)}V^{1-\frac{\alpha}{2}}(t)\|u_{q}\|_{L^{\rho}_{t}L^{p}}\\
        & \quad +2^{(q-j)(1+\frac{\alpha}{\rho})}\int_{0}^{t}c_{q}(\tau)Z'(\tau)
            e^{CV(\tau)}\|u(\tau)\|_{\dot{B}^{s}_{p,r}}\intd \tau.
    \end{split}
\end{equation}
Let $M_{0}\in \mathbb{Z}$ to be fixed hereafter. Decomposing
\begin{equation*}
    u_{q}=\dot{S}_{q-M_{0}}\bar{u}_{q}\circ \psi_{q}^{-1}
    +\sum_{j\geq q-M_{0}}\dot{\Delta}_{j}\bar{u}_{q}\circ \psi_{q}^{-1},
\end{equation*}
we have for all $t\in[0,T]$,
\begin{equation}\label{eq3.7}
    \|u_{q}\|_{L^{\rho}_{t}L^{p}}\leq e^{CV(t)}\big(\|\dot{S}_{q-M_{0}}\bar{u}_{q}\|_{L^{\rho}_{t}L^{p}}
    +\sum_{j\geq q-M_{0}}\|\dot{\Delta}_{j}\bar{u}_{q}\|_{L^{\rho}_{t}L^{p}}\big).
\end{equation}
By Lemma A.1 in \cite{Danchin}, we have
\begin{equation*}
    \|\dot{S}_{q-M_{0}}\bar{u}_{q}\|_{L^{p}}
    \lesssim \|J_{\psi_{q}^{-1}}\|_{L^{\infty}}^{\frac{1}{p}}(2^{-q}\|\nabla J_{\psi_{q}^{-1}}\|_{L^{\infty}}
    \|J_{\psi_{q}}\|_{L^{\infty}}+2^{-M_{0}}\|\nabla \psi_{q}^{-1}\|_{L^{\infty}})\|u_{q}\|_{L^{p}}.
\end{equation*}
This together with Lemma \ref{lem2.3} and Bernstein lemma leads to
\begin{equation}\label{eq3.8}
    \|\dot{S}_{q-M_{0}}\bar{u}_{q}\|_{L^{\rho}_{t}L^{p}}
    \lesssim e^{CV(t)}(e^{CV(t)}-1+2^{-M_{0}})\|u_{q}\|_{L^{\rho}_{t}L^{p}}.
\end{equation}
As $\dot{\Delta}_{j}u_{0,q}=0$ for $|j-q|>1$, from \eqref{eq3.6} we get
\begin{equation}\label{eq3.9}
    \begin{split}
      &\sum_{j\geq q-M_{0}}\nu^{\frac{1}{\rho}}2^{q(s+\frac{\alpha}{\rho})}
        \|\dot{\Delta}_{j}\bar{u}_{q}\|_{L^{\rho}_{t}L^{p}}\\
        \lesssim & (1-e^{-\kappa\nu \rho t2^{q\alpha}})^{\frac{1}{\rho}}2^{qs}\|u_{0,q}\|_{L^{p}}
        +\nu^{-\frac{1}{\rho_{1}'}}2^{M_{0}(1+\alpha)}e^{CV(t)}
            2^{q(s-\frac{\alpha}{\rho_{1}'})}\|f_{q}\|_{L^{\rho_{1}}_{t}L^{p}}\\
        & +\nu^{\frac{1}{\rho}}2^{M_{0}\alpha}e^{CV(t)}V^{1-\frac{\alpha}{2}}(t)
            2^{q(s+\frac{\alpha}{\rho})}\|u_{q}\|_{L^{\rho}_{t}L^{p}}\\
        & +2^{M_{0}(1+\alpha)}\int_{0}^{t}c_{q}(\tau)Z'(\tau)e^{CV(\tau)}\|u(\tau)\|_{\dot{B}^{s}_{p,r}}\intd \tau.
    \end{split}
\end{equation}
Plugging \eqref{eq3.8} and \eqref{eq3.9} into \eqref{eq3.7} yields that
\begin{equation*}
    \begin{split}
      \nu^{\frac{1}{\rho}}2^{q(s+\frac{\alpha}{\rho})}\|u_{q}\|_{L^{\rho}_{t}L^{p}} & \leq
        C(1-e^{-\kappa\nu \rho t2^{q\alpha}})^{\frac{1}{\rho}}2^{qs}\|u_{0,q}\|_{L^{p}}\\
        & \quad +Ce^{CV(t)}\Big(\nu^{-\frac{1}{\rho_{1}'}}2^{M_{0}(1+\alpha)}
            2^{q(s-\frac{\alpha}{\rho_{1}'})}\|f_{q}\|_{L^{\rho_{1}}_{t}L^{p}}\\
        & \quad +\big(2^{-M_{0}}+2^{M_{0}\alpha}V^{1-\frac{\alpha}{2}}(t)\big)\nu^{\frac{1}{\rho}}
            2^{q(s+\frac{\alpha}{\rho})}\|u_{q}\|_{L^{\rho}_{t}L^{p}}\\
        & \quad +2^{M_{0}(1+\alpha)}\int_{0}^{t}c_{q}(\tau)Z'(\tau)\|u(\tau)\|_{\dot{B}^{s}_{p,r}}\intd \tau\Big).
    \end{split}
\end{equation*}
Choose $M_{0}$ to be the unique integer such that $2C2^{-M_{0}}\in(\frac{1}{8},\frac{1}{4}]$ and $T_{1}$ to be
the largest real number such that
\begin{equation*}
    T_{1}\leq T \quad\text{and}\quad CV(T_{1})\leq C_{0} \quad\text{with}\quad
    C_{0}=\min\bigg(\ln 2,\Big(\frac{2^{-M_{0}\alpha}}{8C^{\frac{\alpha}{2}}}\Big)^{\frac{2}{2-\alpha}}\bigg).
\end{equation*}
Thus for $t\in[0,T_{1}]$, there exists a constant $C_{1}$ such that
\begin{equation*}
    \begin{split}
      \nu^{\frac{1}{\rho}}2^{q(s+\frac{\alpha}{\rho})}\|u_{q}\|_{L^{\rho}_{t}L^{p}} & \leq C_{1}
        \Big((1-e^{-\kappa\nu \rho t2^{q\alpha}})^{\frac{1}{\rho}}2^{qs}\|u_{0,q}\|_{L^{p}}\\
        & \quad +\nu^{-\frac{1}{\rho_{1}'}}2^{q(s-\frac{\alpha}{\rho_{1}'})}\|f_{q}\|_{L^{\rho_{1}}_{t}L^{p}}\\
        & \quad +\int_{0}^{t}c_{q}(\tau)Z'(\tau)\|u(\tau)\|_{\dot{B}^{s}_{p,r}}\intd \tau\Big).
    \end{split}
\end{equation*}
Taking $\ell^{r}$ norm yields
\begin{equation}\label{eq3.10}
      \nu^{\frac{1}{\rho}}\|u\|_{\widetilde{L}^{\rho}_{t}\dot{B}^{s+\frac{\alpha}{\rho}}_{p,r}}
      \leq C_{1}\Big(\|u_{0}\|_{\dot{B}^{s}_{p,r}}
      +\nu^{-\frac{1}{\rho_{1}'}}\|f\|_{\widetilde{L}^{\rho_{1}}_{t}\dot{B}^{s-\frac{\alpha}{\rho_{1}'}}_{p,r}}
      +\int_{0}^{t}Z'(\tau)\|u(\tau)\|_{\dot{B}^{s}_{p,r}}\intd \tau\Big).
\end{equation}
Splitting $[0,T]$ into $m$ subintervals like as $[0,T_{1}]$, $[T_{1},T_{2}]$ and so on, such that
\begin{equation*}
    C\int_{T_{k}}^{T_{k+1}}\|\nabla v(t)\|_{L^{\infty}}\intd t\approx C_{0}.
\end{equation*}
Arguing similarly as in deriving \eqref{eq3.10}, we get for all $t\in[T_{k},T_{k+1}]$,
\begin{equation*}
    \begin{split}
      \nu^{\frac{1}{\rho}}\|u\|_{\widetilde{L}^{\rho}_{[T_{k},t]}\dot{B}^{s+\frac{\alpha}{\rho}}_{p,r}}
      & \leq C_{1}\Big(\|u(T_{k})\|_{\dot{B}^{s}_{p,r}}
      +\nu^{-\frac{1}{\rho_{1}'}}\|f\|_{\widetilde{L}^{\rho_{1}}_{[T_{k},t]}\dot{B}^{s-\frac{\alpha}{\rho_{1}'}}_{p,r}}\\
      & \quad +\int_{T_{k}}^{t}Z'(\tau)\|u(\tau)\|_{\dot{B}^{s}_{p,r}}\intd \tau\Big).
    \end{split}
\end{equation*}
By a standard induction argument, it can be shown that
\begin{equation*}
      \nu^{\frac{1}{\rho}}\|u\|_{\widetilde{L}^{\rho}_{t}\dot{B}^{s+\frac{\alpha}{\rho}}_{p,r}}
      \leq C_{1}^{k+1}\Big(\|u_{0}\|_{\dot{B}^{s}_{p,r}}
      +\nu^{-\frac{1}{\rho_{1}'}}\|f\|_{\widetilde{L}^{\rho_{1}}_{t}\dot{B}^{s-\frac{\alpha}{\rho_{1}'}}_{p,r}}
      +\int_{0}^{t}Z'(\tau)\|u(\tau)\|_{\dot{B}^{s}_{p,r}}\intd \tau\Big).
\end{equation*}
Since the number of such subintervals is $m\approx CV(T)C_{0}^{-1}$, one can readily conclude that
up to a change of $C$,
\begin{equation}\label{eq3.11}
    \begin{split}
      \nu^{\frac{1}{\rho}}\|u\|_{\widetilde{L}^{\rho}_{T}\dot{B}^{s+\frac{\alpha}{\rho}}_{p,r}}
      & \leq Ce^{CV(T)}\Big(\|u_{0}\|_{\dot{B}^{s}_{p,r}}
      +\nu^{-\frac{1}{\rho_{1}'}}\|f\|_{\widetilde{L}^{\rho_{1}}_{T}\dot{B}^{s-\frac{\alpha}{\rho_{1}'}}_{p,r}}\\
      & \quad +\int_{0}^{T}Z'(t)\|u(t)\|_{\dot{B}^{s}_{p,r}}\intd t\Big).
    \end{split}
\end{equation}
Of course, the above inequality is valid for all $\rho\in[\rho_{1},\infty]$.
Choosing first $\rho=\infty$ in \eqref{eq3.11} and applying Gronwall lemma leads to
\begin{equation}\label{eq3.12}
    \|u\|_{\widetilde{L}^{\infty}_{T}\dot{B}^{s}_{p,r}}
    \leq Ce^{CZ(T)}\Big(\|u_{0}\|_{\dot{B}^{s}_{p,r}}+\nu^{\frac{1}{\rho_{1}}-1}
    \|f\|_{\widetilde{L}^{\rho_{1}}_{T}\dot{B}^{s-\alpha+\frac{\alpha}{\rho_{1}}}_{p,r}}\Big).
\end{equation}
Now plugging \eqref{eq3.12} into \eqref{eq3.11} yields the desired estimate for general $\rho$.
\end{proof}
\section{Local well-posedness}
In this section, we prove the following result:
\begin{proposition}\label{prop4.1}
Let $u_{0}\in \dot{B}^{\frac{1}{p}}_{p,1}(\mathbb{R})$ with $p\in[1,\infty)$,
then there exists $T>0$ such that the equation \eqref{eq1.2} has a unique solution $u$ such that
\begin{equation*}
    u\in \widetilde{L}^{\infty}_{T}\dot{B}^{\frac{1}{p}}_{p,1}
    \cap L^{1}_{T}\dot{B}^{\frac{1}{p}+1}_{p,1}.
\end{equation*}
Besides for all $\beta\in\mathbb{R}^{+}$,
we have $t^{\beta}u\in\widetilde{L}^{\infty}_{T}\dot{B}^{\frac{1}{p}+\beta}_{p,1}$.
\end{proposition}

\begin{proof}
We prove this proposition by making use of an iterative method.

\vskip0.3cm
\textbf{Step 1}: approximation solution.

Let $u^{0}:=e^{-t\Lambda}u_{0}(x)$ and let $u^{n+1}$ be the solution of the linear equation
\begin{equation*}
    \begin{cases}
        \partial_{t}u^{n+1}+u^{n}\partial_{x}u^{n+1}+\Lambda u^{n+1}=0\\
        u^{n+1}(x,0)=u_{0}(x).
    \end{cases}
\end{equation*}
Obviously $u^{0}\in L^{1}(\mathbb{R}^{+};\dot{B}^{\frac{1}{p}+1}_{p,1})$,
thus according to Theorem \ref{thm1.2}, we have $\forall n\in\mathbb{N}$,
\begin{equation*}
    u^{n}\in \widetilde{L}^{\infty}(\mathbb{R}^{+};\dot{B}^{\frac{1}{p}}_{p,1})
    \cap L^{1}(\mathbb{R}^{+};\dot{B}^{\frac{1}{p}+1}_{p,1}).
\end{equation*}

\vskip0.3cm
\textbf{Step 2}: uniform bounds.

Now we intend to obtain uniform bounds, with respect to the parameter $n$, for some $T>0$ independent of $n$.

By making use of Lemma \ref{lem2.5} and similar arguments as in the proof of Theorem \ref{thm1.2},
for all $T>0$ such that
\begin{equation*}
    \int_{0}^{T}\|u^{n}(\tau)\|_{\dot{B}^{\frac{1}{p}+1}_{p,1}}\intd \tau\leq CC_{0},
\end{equation*}
we have
\begin{equation*}
    \begin{split}
    \|u^{n+1}\|_{\widetilde{L}^{2}_{t}\dot{B}^{\frac{1}{p}+\frac{1}{2}}_{p,1}}
    +\|u^{n+1}\|_{L^{1}_{t}\dot{B}^{\frac{1}{p}+1}_{p,1}}
    &\leq C\sum_{q\in\mathbb{Z}}(1-e^{-\kappa t2^{q}})^{\frac{1}{2}}
        2^{\frac{q}{p}}\|\dot{\Delta}_{q}u_{0}\|_{L^{p}}\\
    &\quad+C\|u^{n}\|_{\widetilde{L}^{2}_{t}\dot{B}^{\frac{1}{p}+\frac{1}{2}}_{p,1}}
    \|u^{n+1}\|_{\widetilde{L}^{2}_{t}\dot{B}^{\frac{1}{p}+\frac{1}{2}}_{p,1}}.
    \end{split}
\end{equation*}
By Lebesgue theorem, there exist $T>0$ and an absolute constant $\varepsilon_{0}>0$ such that
\begin{equation}\label{eq4.1}
    \sum_{q\in\mathbb{Z}}(1-e^{-\kappa T2^{q}})^{\frac{1}{2}}2^{\frac{q}{p}}\|\dot{\Delta}_{q}u_{0}\|_{L^{p}}
    \leq \varepsilon_{0}
\end{equation}
and
\begin{equation}\label{eq4.2}
    \|u^{n+1}\|_{\widetilde{L}^{2}_{T}\dot{B}^{\frac{1}{p}+\frac{1}{2}}_{p,1}}
    +\|u^{n+1}\|_{L^{1}_{T}\dot{B}^{\frac{1}{p}+1}_{p,1}}
    \leq 2\varepsilon_{0}.
\end{equation}
On the other hand, by Theorem \ref{thm1.2}
and the Sobolev embedding $\dot{B}^{\frac{1}{p}}_{p,1}\hookrightarrow L^{\infty}$, we have
\begin{equation*}
    \|u^{n+1}\|_{\widetilde{L}^{\infty}_{T}\dot{B}^{\frac{1}{p}}_{p,1}}
    \leq Ce^{C\int_{0}^{T}\|u^{n}(\tau)\|_{\dot{B}^{\frac{1}{p}+1}_{p,1}}\intd \tau}
    \|u_{0}\|_{\dot{B}^{\frac{1}{p}}_{p,1}}
    \leq C\|u_{0}\|_{\dot{B}^{\frac{1}{p}}_{p,1}}.
\end{equation*}
Combining the above results, we have proved that the sequence $(u^{n})_{n\in \mathbb{N}}$ is uniformly bounded in
$\widetilde{L}^{\infty}_{T}\dot{B}^{\frac{1}{p}}_{p,1}\cap L^{1}_{T}\dot{B}^{\frac{1}{p}+1}_{p,1}$.

\vskip0.3cm
\textbf{Step 3}: strong convergence.

We first prove that $(u^{n})_{n\in \mathbb{N}}$ is a Cauchy sequence
in $\widetilde{L}^{\infty}_{T}\dot{B}^{\frac{1}{p}}_{p,1}$.

Let $(n,m)\in\mathbb{N}^{2}$, $n>m$ and $u^{n,m}:=u^{n}-u^{m}$.
One easily verifies that
\begin{equation*}
    \begin{cases}
        \partial_{t}u^{n+1,m+1}+u^{n}\partial_{x}u^{n+1,m+1}+\Lambda u^{n+1,m+1}=-u^{n,m}\partial_{x}u^{m+1}\\
        u^{n+1,m+1}(x,0)=0.
    \end{cases}
\end{equation*}
According to Theorem \ref{thm1.2}, we have
\begin{equation}\label{eq4.3}
    \|u^{n+1,m+1}\|_{\widetilde{L}^{\infty}_{T}\dot{B}^{\frac{1}{p}}_{p,1}}
    \leq Ce^{C\|u^{n}\|_{L^{1}_{T}\dot{B}^{\frac{1}{p}+1}_{p,1}}}
    \int_{0}^{T}\|u^{n,m}\partial_{x}u^{m+1}(\tau)\|_{\dot{B}^{\frac{1}{p}}_{p,1}}\intd\tau.
\end{equation}
By Proposition \ref{prop2.1} and the embedding $\dot{B}^{\frac{1}{p}}_{p,1}\hookrightarrow L^{\infty}$, we have
\begin{equation*}
    \|u^{n,m}\partial_{x}u^{m+1}\|_{\dot{B}^{\frac{1}{p}}_{p,1}}\lesssim\|u^{n,m}\|_{\dot{B}^{\frac{1}{p}}_{p,1}}
        \|u^{m+1}\|_{\dot{B}^{\frac{1}{p}+1}_{p,1}}.
\end{equation*}
Substituting this into \eqref{eq4.3} yields
\begin{equation*}
    \|u^{n+1,m+1}\|_{\widetilde{L}^{\infty}_{T}\dot{B}^{\frac{1}{p}}_{p,1}}
    \leq C\|u^{n,m}\|_{\widetilde{L}^{\infty}_{T}\dot{B}^{\frac{1}{p}}_{p,1}}
    e^{C\|u^{n}\|_{L^{1}_{T}\dot{B}^{\frac{1}{p}+1}_{p,1}}}
    \int_{0}^{T}\|u^{m+1}(\tau)\|_{\dot{B}^{\frac{1}{p}+1}_{p,1}}\intd\tau.
\end{equation*}
By \eqref{eq4.1}, we can choose $\varepsilon_{0}$ small enough such that
\begin{equation*}
    \|u^{n+1,m+1}\|_{\widetilde{L}^{\infty}_{T}\dot{B}^{\frac{1}{p}}_{p,1}}
    \leq \epsilon\|u^{n,m}\|_{\widetilde{L}^{\infty}_{T}\dot{B}^{\frac{1}{p}}_{p,1}}
\end{equation*}
with $\epsilon<1$.
Now we can get by induction
\begin{equation*}
    \|u^{n+1,m+1}\|_{\widetilde{L}^{\infty}_{T}\dot{B}^{\frac{1}{p}}_{p,1}}
    \leq \epsilon^{m+1}\|u^{n,0}\|_{\widetilde{L}^{\infty}_{T}\dot{B}^{\frac{1}{p}}_{p,1}}
    \leq C\epsilon^{m+1}\|u_{0}\|_{\dot{B}^{\frac{1}{p}}_{p,1}}.
\end{equation*}
This implies that $(u^{n})_{n\in \mathbb{N}}$ is a Cauchy sequence
in $\widetilde{L}^{\infty}_{T}\dot{B}^{\frac{1}{p}}_{p,1}$.
Thus there exists $u\in\widetilde{L}^{\infty}_{T}\dot{B}^{\frac{1}{p}}_{p,1}$ such that $u^{n}$
converges strongly to $u$ in $\widetilde{L}^{\infty}_{T}\dot{B}^{\frac{1}{p}}_{p,1}$.
Fatou lemma and \eqref{eq4.2} ensure that $u\in L^{1}_{T}\dot{B}^{\frac{1}{p}+1}_{p,1}$.
Thus by passing to the limit into the approximation equation, we can get a solution to \eqref{eq1.2} in
$\widetilde{L}^{\infty}_{T}\dot{B}^{\frac{1}{p}}_{p,1}\cap L^{1}_{T}\dot{B}^{\frac{1}{p}+1}_{p,1}$.

\vskip0.3cm
\textbf{Step 4}: uniqueness.

Let $u_{1}$ and $u_{2}$ be two solutions of the equation \eqref{eq1.2} with the same initial data
and belonging to the space
$\widetilde{L}^{\infty}_{T}\dot{B}^{\frac{1}{p}}_{p,1}\cap L^{1}_{T}\dot{B}^{\frac{1}{p}+1}_{p,1}$.
Let $u_{1,2}:=u_{1}-u_{2}$, then we have
\begin{equation*}
    \begin{cases}
        \partial_{t}u_{1,2}+u_{1}\partial_{x}u_{1,2}+\Lambda u_{1,2}=-u_{1,2}\partial_{x}u_{2}\\
        u_{1,2}(x,0)=0.
    \end{cases}
\end{equation*}
By similar arguments as in Step 3, we have
\begin{equation*}
    \|u_{1,2}\|_{\widetilde{L}^{\infty}_{t}\dot{B}^{\frac{1}{p}}_{p,1}}
    \leq Ce^{C\|u_{1}\|_{L^{1}_{t}\dot{B}^{\frac{1}{p}+1}_{p,1}}}
    \int_{0}^{t}\|u_{1,2}\|_{\widetilde{L}^{\infty}_{\tau}\dot{B}^{\frac{1}{p}}_{p,1}}
    \|u_{2}(\tau)\|_{\dot{B}^{\frac{1}{p}+1}_{p,1}}\intd\tau.
\end{equation*}
Gronwall's inequality ensures that $u_{1}=u_{2}$, $\forall t\in[0,T]$.

\vskip0.3cm
\textbf{Step 5}: smoothing effect.

We will prove that for all $\beta\in \mathbb{R}^{+}$, we have
\begin{equation}\label{eq4.4}
    \|t^{\beta}u\|_{\widetilde{L}^{\infty}_{T}\dot{B}^{\frac{1}{p}+\beta}_{p,1}}
    \leq C_{\beta}e^{C(\beta+1)\|u\|_{L^{1}_{T}\dot{B}^{\frac{1}{p}+1}_{p,1}}}
    \|u\|_{\widetilde{L}^{\infty}_{T}\dot{B}^{\frac{1}{p}}_{p,1}}.
\end{equation}
It is obvious that
\begin{equation*}
    \begin{cases}
        \partial_{t}(t^{\beta}u)+u\partial_{x}(t^{\beta}u)+\Lambda (t^{\beta}u)=\beta t^{\beta-1}u\\
        (t^{\beta}u)(x,0)=0.
    \end{cases}
\end{equation*}
When $\beta=1$, by Theorem \ref{thm1.2}, we have
\begin{equation*}
    \|tu(t)\|_{\widetilde{L}^{\infty}_{T}\dot{B}^{\frac{1}{p}+1}_{p,1}}
    \leq Ce^{C\|u\|_{L^{1}_{T}\dot{B}^{\frac{1}{p}+1}_{p,1}}}
    \|u\|_{\widetilde{L}^{\infty}_{T}\dot{B}^{\frac{1}{p}}_{p,1}}.
\end{equation*}
Suppose \eqref{eq4.4} is true for $n$, we will prove it for $n+1$.
Applying Theorem \ref{thm1.2} to the equation of $t^{n+1}u$ yields that
\begin{equation*}
    \begin{split}
      \|t^{n+1}u(t)\|_{\widetilde{L}^{\infty}_{T}\dot{B}^{\frac{1}{p}+n+1}_{p,1}} & \leq C(n+1)
        e^{C\|u\|_{L^{1}_{T}\dot{B}^{\frac{1}{p}+1}_{p,1}}}
        \|t^{n}u\|_{\widetilde{L}^{\infty}_{T}\dot{B}^{\frac{1}{p}+n}_{p,1}} \\
        & \leq C_{n}e^{C(n+2)\|u\|_{L^{1}_{T}\dot{B}^{\frac{1}{p}+1}_{p,1}}}
        \|u\|_{\widetilde{L}^{\infty}_{T}\dot{B}^{\frac{1}{p}}_{p,1}}.
    \end{split}
\end{equation*}
For general $\beta\in \mathbb{R}^{+}$, obviously $[\beta]\leq \beta\leq [\beta]+1$.
Thus by the following interpolation
\begin{equation*}
    \|t^{\beta}u\|_{\widetilde{L}^{\infty}_{T}\dot{B}^{\frac{1}{p}+\beta}_{p,1}}
    \lesssim \|t^{[\beta]}u\|_{\widetilde{L}^{\infty}_{T}\dot{B}^{\frac{1}{p}+[\beta]}_{p,1}}^{[\beta]+1-\beta}
    \|t^{[\beta]+1}u\|_{\widetilde{L}^{\infty}_{T}\dot{B}^{\frac{1}{p}+[\beta]+1}_{p,1}}^{\beta-[\beta]},
\end{equation*}
we can get the estimate for general $\beta\in \mathbb{R}^{+}$.
\end{proof}
\section{Blow-up criterion}
In this section, we prove the following blow-up criterion:
\begin{proposition}\label{prop5.1}
Let $T^{\ast}$ be the maximum local existence time of $u$ in
$\widetilde{L}^{\infty}_{T}\dot{B}^{\frac{1}{p}}_{p,1}\cap
L^{1}_{T}\dot{B}^{\frac{1}{p}+1}_{p,1}$. If $T^{\ast}<\infty$,
then
\begin{equation*}
    \int_{0}^{T^{\ast}}\|\partial_{x}u(t)\|_{L^{\infty}}\intd t=\infty.
\end{equation*}
\end{proposition}

\begin{proof}
Suppose $\int_{0}^{T^{\ast}}\|\partial_{x}u(t)\|_{L^{\infty}}\intd
t$ be finite, then by Theorem \ref{thm1.2}, we have
\begin{equation}\label{eq5.1}
    \forall t\in[0,T^{\ast}),\quad\|u(t)\|_{\dot{B}^{\frac{1}{p}}_{p,1}}
    \leq M_{T^{\ast}}:=Ce^{C\int_{0}^{T^{\ast}}\|\partial_{x}u(t)\|_{L^{\infty}}\intd t}
    \|u_{0}\|_{\dot{B}^{\frac{1}{p}}_{p,1}}<\infty.
\end{equation}
Let $\widetilde{T}>0$ such that
\begin{equation}\label{eq5.2}
    \sum_{q\in\mathbb{Z}}(1-e^{-\kappa \widetilde{T}2^{q}})^{\frac{1}{2}}M_{T^{\ast}}\leq \varepsilon_{0},
\end{equation}
where $\varepsilon_{0}$ is the absolute constant emerged in the
proof of Proposition \ref{prop4.1}. Now \eqref{eq5.1} and
\eqref{eq5.2} imply that
\begin{equation*}
    \forall t\in[0,T^{\ast}),\quad\sum_{q\in\mathbb{Z}}(1-e^{-\kappa \widetilde{T}2^{q}})^{\frac{1}{2}}
    2^{\frac{q}{p}}\|\dot{\Delta}_{q}u(t)\|_{L^{p}}
    \leq \varepsilon_{0}.
\end{equation*}
This together with the local existence theory ensures that, there
exists a solution $\widetilde{u}(t)$ on $[0,\widetilde{T})$ to
\eqref{eq1.2} with the initial datum
$u(T^{\ast}-\widetilde{T}/2)$. By uniqueness,
$\widetilde{u}(t)=u(t+T^{\ast}-\widetilde{T}/2)$ on
$[0,\widetilde{T}/2)$ so that $\widetilde{u}$ extends the solution
$u$ beyond $T^{\ast}$.
\end{proof}
\section{Global well-posedness}
In this section, making use of the method of modulus of continuity \cite{K-N-S},
with help of similar arguments as in \cite{A-H}, we give the proof of the global well-posedness.

Let $T^{\ast}$ be the maximal existence time of the solution $u$ to \eqref{eq1.2} in the space
$\widetilde{L}^{\infty}([0,T^{\ast});\dot{B}^{\frac{1}{p}}_{p,1})
\cap L^{1}_{loc}([0,T^{\ast});\dot{B}^{\frac{1}{p}+1}_{p,1})$.
From Proposition \ref{prop4.1}, there exists $T_{0}>0$ such that
\begin{equation*}
    \forall t\in[0,T_{0}],\quad t\|\partial_{x}u(t)\|_{L^{\infty}}\leq C\|u^{0}\|_{\dot{B}^{\frac{1}{p}}_{p,1}}.
\end{equation*}
Let $\lambda$ be a positive real number that will be fixed later and $T_{1}\in(0,T_{0})$.
We define the set
\begin{equation*}
I:=\{T\in[T_{1},T^{\ast});\forall t\in[T_{1},T],\forall x\neq y\in\mathbb{R}, |u(x,t)-u(y,t)|<\omega_{\lambda}(|x-y|)\},
\end{equation*}
where $\omega:\mathbb{R}^{+}\longrightarrow\mathbb{R}^{+}$ is strictly increasing, concave,
$\omega(0)=0$, $\omega'(0)<+\infty$, $\lim_{\xi\rightarrow 0^{+}}\omega''(\xi)=-\infty$ and
\begin{equation*}
\omega_{\lambda}(|x-y|)=\omega(\lambda|x-y|).
\end{equation*}
The function $\omega$ is a modulus of continuity chosen as in \cite{K-N-S}.

We first prove that $T_{1}$ belongs to $I$ under suitable conditions over $\lambda$.
Let $C_{0}$ be a large positive number such that
\begin{equation}\label{eq6.1}
    2\|u_{0}\|_{L^{\infty}}<\omega(C_{0})<3\|u_{0}\|_{L^{\infty}}.
\end{equation}
Since $\omega$ is strictly increasing, then by maximum principle we have
\begin{equation*}
    \lambda|x-y|\geq C_{0}\Rightarrow |u(x,T_{1})-u(y,T_{1})|\leq 2\|u_{0}\|_{L^{\infty}}
    <\omega_{\lambda}(|x-y|).
\end{equation*}
On the other hand we have from Mean Value Theorem
\begin{equation*}
    |u(x,T_{1})-u(y,T_{1})|\leq|x-y|\|\partial_{x}u(T_{1})\|_{L^{\infty}}.
\end{equation*}
Let $0<\delta_{0}<C_{0}$. Then by the concavity of $\omega$ we have
\begin{equation*}
    \lambda|x-y|\leq \delta_{0}\Rightarrow \omega_{\lambda}(|x-y|)
    \geq \frac{\omega(\delta_{0})}{\delta_{0}}\lambda|x-y|.
\end{equation*}
If we choose $\lambda$ so that
\begin{equation*}
    \lambda>\frac{\delta_{0}}{\omega(\delta_{0})}\|\partial_{x}u(T_{1})\|_{L^{\infty}},
\end{equation*}
then we get
\begin{equation*}
    0<\lambda|x-y|\leq \delta_{0}\Rightarrow |u(x,T_{1})-u(y,T_{1})|<\omega_{\lambda}(|x-y|).
\end{equation*}
Let us now consider the case $\delta_{0}\leq\lambda|x-y|\leq C_{0}$.
By Mean Value Theorem and the increasing property of $\omega$, we can get
\begin{equation*}
    |u(x,T_{1})-u(y,T_{1})|\leq\frac{C_{0}}{\lambda}\|\partial_{x}u(T_{1})\|_{L^{\infty}}\quad\text{and}\quad
    \omega(\delta_{0})\leq \omega_{\lambda}(|x-y|).
\end{equation*}
Choosing $\lambda$ such that
\begin{equation*}
    \lambda>\frac{C_{0}}{\omega(\delta_{0})}\|\partial_{x}u(T_{1})\|_{L^{\infty}},
\end{equation*}
thus we get
\begin{equation*}
    \delta_{0}\leq\lambda|x-y|\leq C_{0}\Rightarrow |u(x,T_{1})-u(y,T_{1})|<\omega_{\lambda}(|x-y|).
\end{equation*}
All the preceding conditions over $\lambda$ can be obtained if we take
\begin{equation}\label{eq6.2}
    \lambda=\frac{\omega^{-1}(3\|u_{0}\|_{L^{\infty}})}{2\|u_{0}\|_{L^{\infty}}}\|\partial_{x}u(T_{1})\|_{L^{\infty}}.
\end{equation}

From the construction, the set $I$ is an interval of the form $[T_{1}, T_{\ast})$.
We have three possibilities which will be discussed separately.

\vskip0.3cm
\textbf{Case 1:}
The first possibility is $T_{\ast}=T^{\ast}$.
In this case we necessarily have $T^{\ast}=\infty$ because the Lipschitz norm of $u$ does not blow up.

\vskip0.3cm
\textbf{Case 2:}
The second possibility is $T_{\ast}\in I$ and we will show that is not possible.

Let $C_{0}$ satisfy \eqref{eq6.1}, then for all $t\in[T_{1},T^{\ast})$, we have
\begin{equation*}
    \lambda|x-y|\geq C_{0}\Rightarrow |u(x,t)-u(y,t)|<\omega_{\lambda}(|x-y|).
\end{equation*}
Since $\partial_{x}u(t)$ belongs to $\mathcal{C}((0,T^{\ast});\dot{B}^{\frac{1}{p}}_{p,1})$,
then for $\varepsilon>0$ there exist $\eta_{0},R>0$ such that $\forall t\in[T_{\ast},T_{\ast}+\eta_{0}]$,
\begin{equation*}
    \|\partial_{x}u(t)\|_{L^{\infty}}\leq\|\partial_{x}u(T_{\ast})\|_{L^{\infty}}
    +\frac{\varepsilon}{2}\quad\text{and}\quad \|\partial_{x}u(T_{\ast})\|_{L^{\infty}(B^{c}_{(0,R)})}
    \leq\frac{\varepsilon}{2},
\end{equation*}
where $B_{(0,R)}$ is the ball of radius $R$ and with center the origin.
Hence for $\lambda|x-y|\leq C_{0}$ and $x$ or $y\in B^{c}_{(0,R+\frac{C_{0}}{\lambda})}$,
we have for $\forall t\in[T_{\ast},T_{\ast}+\eta_{0}]$
\begin{equation*}
    |u(x,t)-u(y,t)|\leq|x-y|\|\partial_{x}u(t)\|_{L^{\infty}(B^{c}_{(0,R)})}\leq\varepsilon|x-y|.
\end{equation*}
On the other hand we have from the concavity of $\omega$
\begin{equation*}
    \lambda|x-y|\leq C_{0}\Rightarrow \frac{\omega(C_{0})}{C_{0}}\lambda|x-y|\leq\omega_{\lambda}(|x-y|).
\end{equation*}
Thus if we take $\varepsilon$ sufficiently small such that
\begin{equation*}
    \varepsilon<\frac{\omega(C_{0})}{C_{0}}\lambda,
\end{equation*}
then we find that
\begin{equation*}
    \lambda|x-y|\leq C_{0};x\;\text{or}\;y\in B^{c}_{(0,R+\frac{C_{0}}{\lambda})}
    \Rightarrow |u(x,t)-u(y,t)|<\omega_{\lambda}(|x-y|).
\end{equation*}

It remains to study the case where $x,y\in B_{(0,R+\frac{C_{0}}{\lambda})}$.
Since $\|\partial^{2}_{x}u(T_{\ast})\|_{L^{\infty}}$ is finite (see Proposition \ref{prop4.1})
then we get for each $x\in\mathbb{R}$
\begin{equation*}
    |\partial_{x}u(x,T_{\ast})|<\lambda\omega'(0).
\end{equation*}
From the continuity of $x\longrightarrow |\partial_{x}u(x,T_{\ast})|$ we obtain
\begin{equation*}
    \|\partial_{x}u(T_{\ast})\|_{L^{\infty}(B_{(0,R+\frac{C_{0}}{\lambda})})}<\lambda\omega'(0).
\end{equation*}
Let $\delta_{1}\ll 1$.
By the continuity in time of the quantity $\|\partial_{x}u(t)\|_{L^{\infty}}$,
there exists $\eta_{1}>0$ such that $\forall t\in[T_{\ast},T_{\ast}+\eta_{1}]$
\begin{equation*}
    \|\partial_{x}u(t)\|_{L^{\infty}(B_{(0,R+\frac{C_{0}}{\lambda})})}<\lambda\frac{\omega(\delta_{1})}{\delta_{1}}.
\end{equation*}
Therefore for $\lambda|x-y|\leq\delta_{1}$ and $x\neq y$ belonging together to $B_{(0,R+\frac{C_{0}}{\lambda})}$,
we have for all $t\in[T_{\ast},T_{\ast}+\eta_{1}]$
\begin{equation*}
    \begin{split}
      |u(x,t)-u(y,t)| & \leq|x-y|\|\partial_{x}u(t)\|_{L^{\infty}(B_{(0,R+\frac{C_{0}}{\lambda})})} \\
        & <\lambda|x-y|\frac{\omega(\delta_{1})}{\delta_{1}}\\
        & \leq \omega_{\lambda}(|x-y|).
    \end{split}
\end{equation*}
Now for the other case since
\begin{equation*}
    \forall x,y\in B_{(0,R+\frac{C_{0}}{\lambda})},\delta_{1}\leq\lambda|x-y|;
    |u(x,T_{\ast})-u(y,T_{\ast})|<\omega_{\lambda}(|x-y|),
\end{equation*}
then we get from a standard compact argument the existence of $\eta_{2}>0$
such that for all $t\in[T_{\ast},T_{\ast}+\eta_{2}]$
\begin{equation*}
    \forall x,y\in B_{(0,R+\frac{C_{0}}{\lambda})},\delta_{1}\leq\lambda|x-y|;
    |u(x,t)-u(y,t)|<\omega_{\lambda}(|x-y|).
\end{equation*}
Taking $\eta=\min(\eta_{0},\eta_{1},\eta_{2})$, we obtain that $T_{\ast}+\eta\in I$
which contradicts the fact that $T_{\ast}$ is maximal.

\vskip0.3cm
\textbf{Case 3:}
The last possibility is that $T_{\ast}$ does not belong to $I$.
By the continuity in time of $u$, there exist $x\neq y$ such that
\begin{equation*}
    u(x,T_{\ast})-u(y,T_{\ast})=\omega_{\lambda}(\xi),\quad\text{with}\quad \xi=|x-y|.
\end{equation*}
We will show that this scenario can not occur and more precisely:
\begin{equation*}
    f'(T_{\ast})<0\quad\text{where}\quad f(t):=u(x,t)-u(y,t).
\end{equation*}
This is impossible since $f(t)\leq f(T_{\ast}),\forall t\in[0,T_{\ast}]$.

The proof is the same as \cite{K-N-S} and for the convenience of the reader we sketch out the proof.
From the regularity of the solution we see that the equation can be defined in the classical manner and
\begin{equation*}
    f'(T_{\ast})=u(y,T_{\ast})\partial_{x}u(y,T_{\ast})-u(x,T_{\ast})\partial_{x}u(x,T_{\ast})
    +\Lambda u(y,T_{\ast})-\Lambda u(x,T_{\ast}).
\end{equation*}
From \cite{K-N-S} we have
\begin{equation*}
    u(y,T_{\ast})\partial_{x}u(y,T_{\ast})-u(x,T_{\ast})\partial_{x}u(x,T_{\ast})
    \leq\omega_{\lambda}(\xi)\omega'_{\lambda}(\xi).
\end{equation*}
Again from \cite{K-N-S}
\begin{equation*}
    \begin{split}
      \Lambda u(y,T_{\ast})-\Lambda u(x,T_{\ast}) & \leq\frac{1}{\pi}\int_{0}^{\frac{\xi}{2}}
        \frac{\omega_{\lambda}(\xi+2\eta)+\omega_{\lambda}(\xi-2\eta)-2\omega_{\lambda}(\xi)}{\eta^{2}}\intd\eta \\
        & +\frac{1}{\pi}\int_{\frac{\xi}{2}}^{\infty}
        \frac{\omega_{\lambda}(2\eta+\xi)-\omega_{\lambda}(2\eta-\xi)-2\omega_{\lambda}(\xi)}{\eta^{2}}\intd\eta\\
        & \leq \lambda J(\lambda\xi),
    \end{split}
\end{equation*}
where
\begin{equation*}
    \begin{split}
    J(\xi)&=\frac{1}{\pi}\int_{0}^{\frac{\xi}{2}}
    \frac{\omega(\xi+2\eta)+\omega(\xi-2\eta)-2\omega(\xi)}{\eta^{2}}\intd\eta\\
    &\quad+\frac{1}{\pi}\int_{\frac{\xi}{2}}^{\infty}
    \frac{\omega(2\eta+\xi)-\omega(2\eta-\xi)-2\omega(\xi)}{\eta^{2}}\intd\eta.
    \end{split}
\end{equation*}
Thus we get
\begin{equation*}
    f'(T_{\ast})\leq\lambda(\omega\omega'+J)(\lambda\xi).
\end{equation*}
Now, we choose the same function as \cite{K-N-S}
\begin{equation*}
    \omega(\xi)=
    \begin{cases}
    \frac{\xi}{1+4\pi\sqrt{\xi_{0}\xi}},& \text{if}\;0\leq\xi\leq \xi_{0};\\
    C_{\xi_{0}}\log\xi,& \text{if}\;\xi\geq \xi_{0},
    \end{cases}
\end{equation*}
here $\xi_{0}$ is sufficiently large number and $C_{\xi_{0}}$ is chosen to provide continuity of $\omega$.
It is shown in \cite{K-N-S},
\begin{equation*}
    \forall \xi\neq 0,\quad\omega(\xi)\omega'(\xi)+J(\xi)<0.
\end{equation*}
Thus we can get that $f'(T_{\ast})<0$.

Combining the above discussion, we conclude that $T^{\ast}=\infty$ and
\begin{equation*}
    \forall t\in[T_{1},\infty),\quad \|\partial_{x}u\|_{L^{\infty}}\leq \lambda \omega'(0)=\lambda.
\end{equation*}
The value of $\lambda$ is given by \eqref{eq6.2}.
\section{Appendix\; -\!- \; Commutator Estimate}
In this appendix, we give the proof of Lemma \ref{lem2.5}.

By Bony's decomposition, we have
\begin{equation}\label{eqA.1}
    \begin{split}
      R_{q} & =(\dot{S}_{q-1}v-v)\cdot\nabla \dot{\Delta}_{q}u-[\dot{\Delta}_{q},v\cdot \nabla]u\\
        & =[\dot{T}_{v^{j}},\dot{\Delta}_{q}]\partial_{j}u
            +\dot{T}_{\partial_{j}\dot{\Delta}_{q}u}v^{j}
            -\dot{\Delta}_{q}\dot{T}_{\partial_{j}u}v^{j}\\
        &\quad +\big\{\partial_{j}\dot{R}(v^{j},\dot{\Delta}_{q}u)
            -\partial_{j}\dot{\Delta}_{q}\dot{R}(v^{j},u)\big\}\\
        &\quad +\big\{\dot{\Delta}_{q}\dot{R}(\Div v,u)
            -\dot{R}(\Div v,\dot{\Delta}_{q}u)\big\}\\
        &\quad +(\dot{S}_{q-1}v-v)\cdot\nabla \dot{\Delta}_{q}u\\
        &=: R_{q}^{1}+R_{q}^{2}+R_{q}^{3}+R_{q}^{4}+R_{q}^{5}+R_{q}^{6}.
    \end{split}
\end{equation}
Above, the summation convention over repeated indices has been used.
The notation $\dot{T}$ stands for homogeneous Bony¡¯s paraproduct which is defined by
\begin{equation*}
    \dot{T}_{f}g:=\sum_{q'\in \mathbb{Z}}\dot{S}_{q'-1}f\dot{\Delta}_{q'}g,
\end{equation*}
and $\dot{R}$ stands for the remainder operator defined by
\begin{equation*}
    \dot{R}(f,g):=\sum_{q'\in \mathbb{Z}}\dot{\Delta}_{q'}f
        (\dot{\Delta}_{q'-1}g+\dot{\Delta}_{q'}g+\dot{\Delta}_{q'+1}g).
\end{equation*}
Note that
\begin{equation}\label{eqA.2}
    \|\dot{\Delta}_{q'}\nabla v\|_{L^{a}}\approx 2^{q'}\|\dot{\Delta}_{q'}v\|_{L^{a}},
    \quad \forall a\in[1,\infty],\; q'\in \mathbb{Z}.
\end{equation}
Now let us estimate each term in \eqref{eqA.1}.

\vskip0.3cm
\textbf{Bounds for $2^{q\sigma}\|R_{q}^{1}\|_{L^{p}}$:}

By \eqref{eq2.1} and the definition of $\dot{\Delta}_{q}$, we have
\begin{equation}\label{eqA.3}
    \begin{split}
      R_{q}^{1} & =\sum_{|q'-q|\leq 4}[\dot{S}_{q'-1}v^{j},\dot{\Delta}_{q}]\partial_{j}\dot{\Delta}_{q'}u\\
        & =\sum_{|q'-q|\leq 4}\int_{\mathbb{R}^{N}} h(y)
            \big[\dot{S}_{q'-1}v^{j}(x)-\dot{S}_{q'-1}v^{j}(x-2^{-q}y)\big]
            \partial_{j}\dot{\Delta}_{q'}u(x-2^{-q}y)\intd y.
    \end{split}
\end{equation}
Applying Mean Value Theorem and Young's inequality to \eqref{eqA.3} yields
\begin{equation}\label{eqA.4}
      2^{q\sigma}\|R_{q}^{1}\|_{L^{p}}\leq C\sum_{|q'-q|\leq 4}\|\dot{S}_{q'-1}\nabla v\|_{L^{\infty}}
            2^{q'\sigma}\|\dot{\Delta}_{q'}u\|_{L^{p}}.
\end{equation}

\vskip0.3cm
\textbf{Bounds for $2^{q\sigma}\|R_{q}^{2}\|_{L^{p}}$:}

According to \eqref{eq2.1}, we have
\begin{equation*}
    R_{q}^{2}=\sum_{q'\geq q-3}\dot{S}_{q'-1}\partial_{j}\dot{\Delta}_{q}u
            \dot{\Delta}_{q'}v^{j}.
\end{equation*}
By \eqref{eqA.2}, we can get
\begin{equation}\label{eqA.5}
    \begin{split}
      2^{q\sigma}\|R_{q}^{2}\|_{L^{p}} & \leq C\sum_{q'\geq q-3}2^{q\sigma}
            \|\dot{\Delta}_{q'}v^{j}\|_{L^{\infty}}
            \|\dot{S}_{q'-1}\partial_{j}\dot{\Delta}_{q}u\|_{L^{p}} \\
        & \leq C\sum_{q'\geq q-3}2^{q-q'}\|\dot{\Delta}_{q'}\nabla v\|_{L^{\infty}}
            2^{q\sigma}\|\dot{\Delta}_{q}u\|_{L^{p}}.
    \end{split}
\end{equation}

\vskip0.3cm
\textbf{Bounds for $2^{q\sigma}\|R_{q}^{3}\|_{L^{p}}$:}

Again from \eqref{eq2.1}, we have
\begin{equation}\label{eqA.6}
    R_{q}^{3}=-\sum_{|q'-q|\leq 4}\dot{\Delta}_{q}(\dot{S}_{q'-1}\partial_{j}u\dot{\Delta}_{q'}v^{j})
    =-\sum_{\substack{|q'-q|\leq 4\\q''\leq q'-2}}
    \dot{\Delta}_{q}(\dot{\Delta}_{q''}\partial_{j}u\dot{\Delta}_{q'}v^{j}).
\end{equation}
Therefore, denoting $\frac{1}{p_{2}}=\frac{1}{p}-\frac{1}{p_{1}}$
and taking advantage of \eqref{eqA.2}, we can obtain
\begin{equation}\label{eqA.7}
    \begin{split}
      2^{q\sigma}\|R_{q}^{3}\|_{L^{p}} & \leq C\sum_{\substack{|q'-q|\leq 4\\q''\leq q'-2}}2^{q\sigma}
      \|\dot{\Delta}_{q'}v^{j}\|_{L^{p_{1}}}\|\dot{\Delta}_{q''}\partial_{j}u\|_{L^{p_{2}}} \\
        & \leq C\sum_{\substack{|q'-q|\leq 4\\q''\leq q'-2}}2^{(q-q'')(\sigma-1-\frac{N}{p_{1}})}
      2^{q'\frac{N}{p_{1}}}\|\dot{\Delta}_{q'}\nabla v\|_{L^{p_{1}}}2^{q''\sigma}\|\dot{\Delta}_{q''}u\|_{L^{p}}.
    \end{split}
\end{equation}
Note that, starting from the first equality of \eqref{eqA.6}, one can alternately get
\begin{equation}\label{eqA.8}
    \begin{split}
      2^{q\sigma}\|R_{q}^{3}\|_{L^{p}} & \leq C\sum_{|q'-q|\leq 4}2^{q\sigma}
      \|\dot{\Delta}_{q'}v^{j}\|_{L^{p_{1}}}\|\dot{S}_{q'-1}\partial_{j}u\|_{L^{p_{2}}} \\
        & \leq C\sum_{|q'-q|\leq 4}2^{q'(\sigma-1)}\|\dot{\Delta}_{q'}\nabla v\|_{L^{p_{1}}}
            \|\dot{S}_{q'-1}\partial_{j}u\|_{L^{p_{2}}}.
    \end{split}
\end{equation}

\vskip0.3cm
\textbf{Bounds for $2^{q\sigma}\|R_{q}^{4}\|_{L^{p}}$:}

\begin{equation*}
    \begin{split}
    R_{q}^{4}&=\sum_{\substack{|q'-q|\leq 2\\|q''-q'|\leq 1}}\partial_{j}(\dot{\Delta}_{q'}v^{j}
        \dot{\Delta}_{q}\dot{\Delta}_{q''}u)
        -\sum_{\substack{q'\geq q-3\\|q''-q'|\leq 1}}\partial_{j}\dot{\Delta}_{q}(\dot{\Delta}_{q'}v^{j}
        \dot{\Delta}_{q''}u)\\
    &=: R_{q}^{4,1}+R_{q}^{4,2}.
    \end{split}
\end{equation*}
By \eqref{eqA.2}, we can get
\begin{equation}\label{eqA.9}
    \begin{split}
        2^{q\sigma}\|R_{q}^{4,1}\|_{L^{p}}&\leq C\sum_{\substack{|q'-q|\leq 2\\|q''-q'|\leq 1}}
            \|\dot{\Delta}_{q'}\nabla v\|_{L^{\infty}}2^{q'\sigma}\|\dot{\Delta}_{q''}u\|_{L^{p}}\\
        &\leq C\sum_{\substack{|q'-q|\leq 2\\|q''-q'|\leq 1}}2^{q'\frac{N}{p_{1}}}
            \|\dot{\Delta}_{q'}\nabla v\|_{L^{p_{1}}}2^{q'\sigma}\|\dot{\Delta}_{q''}u\|_{L^{p}}.
    \end{split}
\end{equation}
For $R_{q}^{4,2}$, we proceed differently according to the value
of $\frac{1}{p}+\frac{1}{p_{1}}$. If
$\frac{1}{p}+\frac{1}{p_{1}}\leq 1$, we denote
$\frac{1}{p_{3}}:=\frac{1}{p}+\frac{1}{p_{1}}$ and have
\begin{equation}\label{eqA.10}
    \begin{split}
      2^{q\sigma}\|R_{q}^{4,2}\|_{L^{p}} & \leq C\sum_{\substack{q'\geq q-3\\|q''-q'|\leq 1}}
        2^{q(1+\sigma)}2^{q(\frac{N}{p_{3}}-\frac{N}{p})}
      \|\dot{\Delta}_{q'}v\dot{\Delta}_{q''}u\|_{L^{p_{3}}}\\
        & \leq C\sum_{\substack{q'\geq q-3\\|q''-q'|\leq 1}}2^{q(1+\sigma)}2^{q\frac{N}{p_{1}}}
      \|\dot{\Delta}_{q'}v\|_{L^{p_{1}}}\|\dot{\Delta}_{q''}u\|_{L^{p}}\\
        & \leq C\sum_{\substack{q'\geq q-3\\|q''-q'|\leq 1}}2^{(q-q')(1+\sigma+\frac{N}{p_{1}})}2^{q'\frac{N}{p_{1}}}
      \|\dot{\Delta}_{q'}\nabla v\|_{L^{p_{1}}}2^{q'\sigma}\|\dot{\Delta}_{q''}u\|_{L^{p}}.
    \end{split}
\end{equation}
If $\frac{1}{p}+\frac{1}{p_{1}}>1$, taking $p_{1}=p'$ in the above computations yields
\begin{equation}\label{eqA.11}
    \begin{split}
      2^{q\sigma}\|R_{q}^{4,2}\|_{L^{p}} & \leq C\sum_{\substack{q'\geq q-3\\|q''-q'|\leq 1}}
        2^{q(1+\sigma)}2^{q\frac{N}{p'}}
      \|\dot{\Delta}_{q'}v\dot{\Delta}_{q''}u\|_{L^{1}}\\
        & \leq C\sum_{\substack{q'\geq q-3\\|q''-q'|\leq 1}}2^{q(1+\sigma)}2^{q\frac{N}{p'}}
      \|\dot{\Delta}_{q'}v\|_{L^{p'}}\|\dot{\Delta}_{q''}u\|_{L^{p}}\\
        & \leq C\sum_{\substack{q'\geq q-3\\|q''-q'|\leq 1}}2^{(q-q')(1+\sigma+\frac{N}{p'})}2^{q'\frac{N}{p_{1}}}
      \|\dot{\Delta}_{q'}\nabla v\|_{L^{p_{1}}}2^{q'\sigma}\|\dot{\Delta}_{q''}u\|_{L^{p}}.
    \end{split}
\end{equation}
Putting \eqref{eqA.9}, \eqref{eqA.10} and \eqref{eqA.11} together, we obtain
\begin{equation}\label{eqA.12}
    \begin{split}
      2^{q\sigma}\|R_{q}^{4}\|_{L^{p}}&\leq C\sum_{\substack{q'\geq q-3\\|q''-q'|\leq 1}}
      2^{(q-q')\big(1+\sigma+N\min(\frac{1}{p_{1}}+\frac{1}{p'})\big)}\\
      &\qquad \times 2^{q'\frac{N}{p_{1}}}
      \|\dot{\Delta}_{q'}\nabla v\|_{L^{p_{1}}}2^{q'\sigma}\|\dot{\Delta}_{q''}u\|_{L^{p}}.
    \end{split}
\end{equation}

\vskip0.3cm
\textbf{Bounds for $2^{q\sigma}\|R_{q}^{5}\|_{L^{p}}$:}

Similar computations yield
\begin{equation}\label{eqA.13}
    \begin{split}
      2^{q\sigma}\|R_{q}^{5}\|_{L^{p}}&\leq C\sum_{\substack{q'\geq q-3\\|q''-q'|\leq 1}}
      2^{(q-q')\big(\sigma+N\min(\frac{1}{p_{1}}+\frac{1}{p'})\big)}\\
      &\qquad \times 2^{q'\frac{N}{p_{1}}}
      \|\dot{\Delta}_{q'}\Div v\|_{L^{p_{1}}}2^{q'\sigma}\|\dot{\Delta}_{q''}u\|_{L^{p}}.
    \end{split}
\end{equation}

\vskip0.3cm
\textbf{Bounds for $2^{q\sigma}\|R_{q}^{6}\|_{L^{p}}$:}

\begin{equation*}
    R_{q}^{6}=-\sum_{q'\geq q-1}\dot{\Delta}_{q'}v\cdot\nabla \dot{\Delta}_{q}u,
\end{equation*}
thus by Bernstein lemma, we have
\begin{equation}\label{eqA.14}
      2^{q\sigma}\|R_{q}^{6}\|_{L^{p}}\leq C\sum_{q'\geq q-1}2^{q-q'}\|\dot{\Delta}_{q'}\nabla v\|_{L^{\infty}}
        2^{q\sigma}\|\dot{\Delta}_{q}u\|_{L^{p}}.
\end{equation}

Combining inequalities \eqref{eqA.4}, \eqref{eqA.5}, \eqref{eqA.7} or \eqref{eqA.8},
\eqref{eqA.12}, \eqref{eqA.13} and \eqref{eqA.14},
we end up with the desired estimate for $R_{q}$.

Straightforward modifications in the estimates for $R_{q}^{3}$,
$R_{q}^{4}$ and $R_{q}^{5}$ leads to the desired estimate in the
special case where $u=v$.

\textbf{Acknowledgments:} The authors would like to thank Prof.
P.Constantin  for helpful comments and suggestions. The authors
also thank Prof. H.Dong and J.Wu for kindly informing us the
recent paper \cite{Dong2}. The authors were partly supported by
the NSF of China (No.10725102).

\end{document}